\documentclass[reqno]{amsart}
\usepackage{graphicx}
\usepackage{verbatim}
\usepackage{textcomp}
\usepackage{amssymb}
\usepackage{cite}
\usepackage{amsmath}
\usepackage{latexsym}
\usepackage{amscd}
\usepackage{amsthm}
\usepackage{mathrsfs}
\usepackage{xypic}
\usepackage{bm}
\usepackage{url}
\usepackage[linktocpage=true,colorlinks,citecolor=magenta,linkcolor=blue,urlcolor=magenta]{hyperref}

\newtheorem{thm}{Theorem}[section]
\newtheorem{corr}[thm]{Corollary}
\newtheorem{lem}[thm]{Lemma}
\newtheorem{prop}[thm]{Proposition}

\theoremstyle{definition}

\newtheorem*{ack}{Acknowledgment}
\theoremstyle{remark}
\newtheorem{rem}{Remark}[section]
\numberwithin{equation}{section}
\begin{document}
\title[Integral pinched shrinking $\rho$-Einstein solitons]
{Integral pinched gradient shrinking $\rho$-Einstein solitons}
\author{Guangyue Huang}
\address{Department of Mathematics, Henan Normal
University, Xinxiang 453007, P.R. China}
\email{\href{mailto:hgy@henannu.edu.cn}{hgy@henannu.edu.cn}}

\thanks{The first author is supported by NSFC(No. 11371018,
11671121).}

\maketitle

\begin{abstract} The gradient shrinking $\rho$-Einstein soliton is a triple $(M^n,g,f)$ such that
$$
R_{ij}+f_{ij}=(\rho R+\lambda) g_{ij},
$$
where $(M^n,g)$ is a Riemannian manifold, $\lambda>0, \rho\in\mathbb{R}\setminus\{0\}$ and $f$ is the potential function on $M^n$. In this paper, using algebraic curvature estimates and the Yamabe-Sobolev inequality,
we prove some integral pinching rigidity results for compact gradient shrinking $\rho$-Einstein solitons.
\end{abstract}

\section{Introduction}

Let $(M^n,g)$ be an $n$-dimensional Riemannian manifold. If there exist two
real constants $\rho\in\mathbb{R}\setminus\{0\},\lambda\in\mathbb{R}$ such that
\begin{equation}\label{Int2}
R_{ij}+f_{ij}=(\rho R+\lambda)g_{ij},
\end{equation}
we call $(M^n,g)$ a gradient $\rho$-Einstein soliton, where $R_{ij}$
is the Ricci curvature and $R$ is the scalar curvature of the metric
$g$, respectively. The gradient $\rho$-Einstein solitons give rise to the self-similar solutions to a perturbed version of the Ricci flow, the so-called Ricci-Bourguignon flows
\begin{equation}\label{rbflow}
  \frac{\partial}{\partial t}g=-2({\rm Ric}-\rho Rg).
\end{equation}
The short-time existence of the flow \eqref{rbflow} was proved in \cite{RBflow} for every $\rho<\frac 1{2(n-1)}$. We refer the readers to \cite{RBflow} for more results on this flow \eqref{rbflow} and \cite{Catino2015,Catino2016N} for results on gradient $\rho$-Einstein solitons. The gradient $\rho$-Einstein soliton $(M^n,g)$ satisfying \eqref{Int2} is called shrinking if $\lambda>0$, steady if $\lambda=0$ and expanding if $\lambda<0$. The gradient $\rho$-Einstein
soliton is a special case of the Ricci almost soliton which is
defined by
\begin{equation}\label{1-addInt2}
R_{ij}+f_{ij}=\tilde{\lambda}\,g_{ij},
\end{equation}
where $\tilde{\lambda}$ is a smooth function on $M^n$, see
\cite{Pigola2011,Bar1}. For more research on equation
\eqref{Int2}, see
\cite{ChenCheng2,Kim03,CaoChen1,Catino2011,Wei2011,HW2013} and the
references therein.

The equation \eqref{Int2} and the flow \eqref{rbflow} are of special interesting for special values of $\rho$: if $\rho=1/n$, $Ric-\rho R$ is just the trace-less Ricci tensor; if $\rho=1/2$, $Ric-\rho R$ is the Einstein tensor; if $\rho=\frac 1{2(n-1)}$, $Ric-\rho R$ is the Schouten tensor. In particular, if $\rho=0$, the equation \eqref{Int2} is exactly the gradient Ricci soliton equation and the flow \eqref{rbflow} is exactly the Hamilton's Ricci flow. As is well known, the gradient Ricci solitons play an important role in Ricci flow as
they correspond to the self-similar solutions and often arise as
singularity models.

In \cite{Catino2016}, Catino proved that for $4\leq n\leq6$, any
compact gradient shrinking Ricci soliton satisfying some
$L^{\frac{n}{2}}$-pinching condition is isometric to a quotient of the round sphere. The proof in \cite{Catino2016}
was inspired by the rigidity of Einstein metrics satisfying a $L^{\frac{n}{2}}$-pinching condition (see \cite{Hebey1996}), and relies on some sharp curvature estimates and the Yamabe-Sobolev
inequality. The main goal of this paper is to classify the gradient $\rho$-Einstein
soliton satisfying some similar $L^{\frac{n}{2}}$-pinching condition and to generalize the results proved
by Catino in \cite{Catino2016}.

Given an $n$-dimensional Riemannian manifold $(M^n,g)$, we denote by $W, \mathrm{Ric},\mathring{\mathrm{Ric}}$ and $R$ the Weyl tensor, Ricci tensor, trace-less part of Ricci tensor and the scalar curvature, respectively. The volume of $M$ with respect to the metric $g$ is denoted by $V(M)$. The Yamabe invariant $Y(M,[g])$ associated to $(M^n,g)$ is defined by
\begin{equation}\label{Int3}\aligned
Y(M,[g])=&\inf\limits_{\tilde{g}\in[g]
}\frac{\int_M\tilde{R}\,dv_{\tilde{g}}}{(\int_M\,dv_{\tilde{g}})^{\frac{n-2}{n}}}\\
=&\frac{4(n-1)}{n-2}\inf\limits_{u\in W^{1,2}(M)}\frac{\int_M|\nabla
u|^2\,dv_g+\frac{n-2}{4(n-1)}\int_M
Ru^2\,dv_g}{(\int_M|u|^{\frac{2n}{n-2}}\,dv_g)^{\frac{n-2}{n}}},
\endaligned\end{equation}
 where $[g]$ is the conformal class of the metric $g$. If $M$ is compact, $Y(M,[g])$ is positive if and only if there exists a conformal metric in $[g]$ with everywhere positive scalar curvature.
As we will show in \S \ref{sec:2}, any compact gradient shrinking $\rho$-Einstein soliton with
$\rho<\frac{1}{2(n-1)}$ has positive scalar curvature, then it has positive Yamabe invariant $Y(M,[g])>0$ as well.

Our first theorem is as follows.
\begin{thm}\label{main-thm1}
Let $(M^n,g)$ be an $n$-dimensional ($4\leq n\leq 6$)  compact gradient shrinking $\rho$-Einstein soliton
satisfying \eqref{Int2} with $\rho<0$. If
\begin{equation}\label{main-ineq1}\aligned
\left(\int\limits_M\Big|W+\frac{\sqrt{2}}{\sqrt{n}(n-2)} \mathring{{\rm
Ric}} \mathbin{\bigcirc\mkern-15mu\wedge}
g\Big|^{\frac{n}{2}}dv_g\right)^{\frac{2}{n}}
+&\sqrt{\frac{(n-4)^{2}(n-1)}{8(n-2)}}\lambda V(M)^{\frac{2}{n}}\\
\leq&~\sqrt{\frac{n-2}{32(n-1)}}Y(M,[g]),
\endaligned\end{equation}
then $M^n$ is isometric to a quotient of the round sphere $\mathbb{S}^n$.
\end{thm}
\begin{rem}
By integrating the equation \eqref{lem1} in \S \ref{sec:2},
\begin{equation*}
  \lambda V(M)^{\frac 2n}=\frac {1-n\rho}nV(M)^{\frac{2-n}n}\int\limits_MR~dv_g\geq \frac {1-n\rho}nY(M,[g])\geq \frac 1nY(M,[g]).
\end{equation*}
Thus, in the case $n\geq 7$, the pinching condition \eqref{main-ineq1} doesn't hold.
\end{rem}

In particular, in the case $n=4$, we have
\begin{corr}\label{1Th1}
Let $(M^4,g)$ be a compact gradient shrinking $\rho$-Einstein soliton
satisfying \eqref{Int2} with $\rho<0$. If
\begin{equation}\label{1th-1}\aligned
\int\limits_M\Big(|W|^2+|\mathring{\rm
Ric}|^2\Big)dv_g<\frac{1}{48}Y(M,[g])^2,
\endaligned\end{equation}
then $M^4$ is isometric to a quotient of the round sphere $\mathbb{S}^4$.
\end{corr}

\begin{rem} By the well-known Chern-Gauss-Bonnet formula:
$$\int\limits_M\Big(|W|^2-2|\mathring{{\rm Ric}}|^2+\frac{1}{6}R^2\Big)dv_g=31\pi^2\chi(M^4),$$
 where $\chi(M)$ is the
Euler-Poincar\'e characteristic of $M^4$, and following the proof of Catino in \cite[\S 4]{Catino2016}, we can also
obtain that, for compact gradient shrinking  $\rho$-Einstein soliton $(M^4,g)$ with
$\rho<0$, if
$$
\int\limits_M|W|^{2}dv_g+\frac{5}{4}\int\limits_M|\mathring{{\rm
Ric}}|^2dv_g~<~\frac{1}{48}\int\limits_MR^{2}dv_g,
$$
then $M^4$ is isometric to a quotient of the  round sphere $\mathbb{S}^4$.
\end{rem}

\begin{rem} Our Theorem \ref{main-thm1} combining with Theorem 1.4 of Catino in \cite{Catino2016} shows
that any compact gradient shrinking $\rho$-Einstein soliton $(M^n,g)$, $4\leq n\leq 6$, with $\rho\leq0$ satisfying \eqref{main-ineq1}  must be isometric to a quotient of the round sphere.
\end{rem}

For the gradient shrinking $\rho$-Einstein soliton $(M^n,g)$ with $\rho>0$, we have the following theorem.
\begin{thm}\label{1Th2}
Let $(M^n,g)$ be an $n$-dimensional ($5\leq n\leq 6$) compact gradient shrinking $\rho$-Einstein soliton satisfying \eqref{Int2}
with positive Ricci curvature. If
\begin{equation}\label{2th-1}\aligned
0<\rho\leq\frac{2(n-2)(n-3)-\sqrt{4(n-2)^2(n-3)^2-n(n-4)^3}}{2n(n-1)(n-4)},
\endaligned\end{equation}
and
\begin{equation}\label{2th-2}\aligned
\Big[\int\limits_M\Big|W+\frac{\sqrt{2}}{\sqrt{n}(n-2)} \mathring{{\rm
Ric}} \mathbin{\bigcirc\mkern-15mu\wedge} g\Big|^{\frac{n}{2}}dv_g\Big]^{\frac{2}{n}}
+&\sqrt{\frac{(n-4)^{2}(n-1)}{8(n-2)}}\lambda V(M)^{\frac{2}{n}}\\
\leq&~\sqrt{\frac{n-2}{32(n-1)}}Y(M,[g]),
\endaligned\end{equation}
and at least one of two inequalities is strict, then
$M^n$ is isometric to a quotient of the  round sphere $\mathbb{S}^n$.
\end{thm}
\begin{rem}
It's easy to check that for $n>4$, if $\rho$ satisfies the condition \eqref{2th-1}, then $0<\rho<\frac 1{2(n-1)}$.
\end{rem}

We observe that the result in Theorem \ref{1Th2} also holds for complete (possibly noncompact) gradient shrinking $\rho$-Einstein soliton with $\rho>0$ and bounded curvature, nonnegative sectional curvature and positive Ricci curvature, since these conditions force the manifold to be compact, as the following theorem says.

\begin{thm}\label{1Th3}
Let $(M^n,g)$ be a complete gradient shrinking $\rho$-Einstein soliton with $0<\rho<\frac 1{2(n-1)}$, bounded curvature, nonnegative sectional curvature and positive Ricci curvature. Then $(M^n,g)$ must be compact.
\end{thm}

Munteanu-Wang \cite{Munt-Wang} already shown that any gradient shrinking Ricci soliton with nonnegative sectional curvature and positive Ricci curvature is compact. Thus Theorem \ref{1Th3} is a generalization of their result to the gradient shrinking $\rho$-Einstein soliton. We remark that the assumption that the curvature is bounded is not required in Munteanu-Wang's theorem. Here this assumption is used to estimate the growth of the potential function, as shown in \cite{Catino2015}.

\begin{ack}
The author would like to express their thanks to Dr. Yong Wei for his helpful discussions and valuable suggestions which make the paper more readable.
\end{ack}

\section{Preliminaries}\label{sec:2}
In this section, we collect some fundamental identities and properties for the gradient $\rho$-Einstein solitons.
\begin{lem}[\cite{Catino2015,Catino2016N}]\label{Lemma21}
Let $(M^n,g)$ be a gradient $\rho$-Einstein soliton satisfying
\eqref{Int2}. Then we have
\begin{align}
&\Delta f=-(1-n\rho) R+n\lambda,\label{lem1}\\
&(1-2(n-1)\rho)R_{,i}=2R_{ij}f^j,\label{lem3}\\
&(1-2(n-1)\rho)\Delta R=\nabla R\nabla f-2|\mathring{{\rm Ric}}|^2-\frac{2}{n}((1-n\rho)R-n\lambda) R,\label{lem4}
\end{align}
where $\mathring{{\rm Ric}}$ is the trace free part of Ricci curvature and $f^j=g^{jk}f_k$.
\end{lem}

From the Lemma \ref{Lemma21} and using the maximum principle, we can obtain the following
\begin{prop}[\cite{Catino2016N}]\label{Prop21} Let $(M^n,g)$ be a compact gradient shrinking $\rho$-Einstein soliton satisfying
\eqref{Int2} with $\rho<\frac{1}{2(n-1)}$, then $(M^n,g)$ has positive scalar curvature $R>0$.
\end{prop}
\proof We include the proof here for convenience of readers. If $R$ is a constant, then \eqref{lem1} shows that
\begin{equation*}
\aligned
0=\int\limits_M\Delta f=\left(n\lambda-(1-n\rho) R\right)V(M)
\endaligned\end{equation*}
which implies that $R={n\lambda}/{(1-n\rho)}>0.$ If $R$ is not a constant, then from \eqref{lem4},
we have
\begin{equation}\label{Int7}\aligned
0\geq~\left((1-n\rho)R_{\min}-n\lambda\right)R_{\min}.
\endaligned\end{equation}
Integrating \eqref{lem1} gives
\begin{equation*}
(1-n\rho)R_{\min}-n\lambda<\frac 1{V(M)}\int_M((1-n\rho)R-n\lambda) = -\frac 1{V(M)}\int_M\Delta f=0.
\end{equation*}
Hence, from \eqref{Int7} we have $R_{\min}\geq0.$ On the other hand, from \eqref{lem4} gives
\begin{equation*}
  (1-2(n-1)\rho)\Delta R-\nabla f\nabla R-2\lambda R\leq 0.
\end{equation*}
By the strong minimum principle, $R$ can not achieve its non-positive minimum in $M^n$
unless $R$ is constant. Hence, we obtain $R_{\min}\neq0$ and $R\geq R_{\min}>0.$
\endproof

\begin{lem}\label{Lemma22}
Let $(M^n,g)$ be a gradient $\rho$-Einstein soliton satisfying
\eqref{Int2}. Then we have
\begin{equation}\label{2lemma-1}\aligned
\frac{1}{2}\Delta_f|\mathring{{\rm Ric}}|^2=&|\nabla\mathring{{\rm Ric}}|^2+2\lambda|\mathring{{\rm Ric}}|^2-2W_{ijkl}\mathring{R}_{ik}\mathring{R}_{jl}\\
&+\frac{4}{n-2}\mathring{R}_{ij}\mathring{R}_{jk}\mathring{R}_{ki}+2\left(\rho-\frac{n-2}{n(n-1)}\right)R|\mathring{{\rm Ric}}|^2+(n-2)\rho R_{,ij}\mathring{R}_{ij}.
\endaligned\end{equation}
\end{lem}
\proof It has been shown in \cite[Lemm 3.3]{{Pigola2011}} that for
a gradient Ricci almost soltion of the form \eqref{1-addInt2}, we
have
\begin{equation*}
\aligned
\Delta_fR_{ik}=&(\Delta\tilde{\lambda})g_{ik}+(n-2)\tilde{\lambda}_{ik}+2\tilde{\lambda}R_{ik}-\frac{2}{n-2}\left(| {\rm Ric}|^2-\frac{R^2}{n-1}\right)g_{ik}\\
&-\frac{2n}{(n-1)(n-2)}RR_{ik}+\frac{4}{n-2}R_{il}R_{kl}-2W_{ijkl}R_{jl}.
\endaligned\end{equation*}
By $\mathring{R}_{ij}=R_{ij}-{R}/{n}g_{ij}$ and
$$R_{ij}R_{jk}R_{ki}=\mathring{R}_{ij}\mathring{R}_{jk}\mathring{R}_{ki}+\frac{3}{n}R|{\rm Ric}|^2-\frac{2}{n^2}R^3,$$
a direction computation shows
\begin{equation*}
\aligned
\frac{1}{2}\Delta_f|\mathring{{\rm Ric}}|^2=&|\nabla\mathring{{\rm Ric}}|^2+2\Big[\tilde{\lambda}-\frac{n-2}{n(n-1)}R\Big]|\mathring{{\rm Ric}}|^2-2W_{ijkl}\mathring{R}_{ik}\mathring{R}_{jl}\\
&+\frac{4}{n-2}\mathring{R}_{ij}\mathring{R}_{jk}\mathring{R}_{ki}+(n-2)\tilde{\lambda}_{ij}\mathring{R}_{ij}.
\endaligned\end{equation*}
Replacing $\tilde{\lambda}$ with $\rho R+\lambda$ yields the desired estimate \eqref{2lemma-1}.
\endproof

We conclude this section with two algebraic curvature estimates involving the Weyl tensor and trace-less Ricci tensor.
\begin{lem}[\cite{Bour,Catino2016}]\label{weyl-est}
On every $n$-dimensional Riemannian manifold $(M^n,g)$, we have
\begin{equation*}
\Big|-W_{ijkl}\mathring{R}_{ik}\mathring{R}_{jl}+\frac{2}{n-2}\mathring{R}_{ij}\mathring{R}_{jk}\mathring{R}_{ki}\Big|
\leq\sqrt{\frac{n-2}{2(n-1)}}\left(|W|^2+\frac{8}{n(n-2)}|\mathring{{\rm
Ric}}|^2\right)^{1/2}|\,\mathring{{\rm Ric}}|^2.
\end{equation*}
\end{lem}

\begin{lem}[\cite{Catino2016,FuXiao2015}]\label{Lemma23}
There exists a positive constant $C(n)$ such that on every $n$-dimensional Riemannian manifold $(M^n,g)$,
the following estimate holds
\begin{equation}\label{1Int10}
2W_{ijkl}W_{ipkq}W_{pjql}+\frac{1}{2}W_{ijkl}W_{klpq}W_{pqij}\leq
C(n)|W|^3,
\end{equation} where
\begin{equation*}
C(n)=
\begin{cases}\frac{\sqrt{6}}{4},\ \ \quad n=4\\
\frac{4\sqrt{10}}{15},  \quad n=5\\
\frac{\sqrt{70}}{2\sqrt{3}}, \ \quad n=6\\
\frac{5}{2}, \ \ \ \ \quad n\geq 7.
\end{cases}
\end{equation*}
\end{lem}

\proof
For the proof,  we refer to \cite[Lemma 2.3]{Catino2016} and \cite[Lemma 2.1]{FuXiao2015}. In \cite{FuXiao2015}, the constant $C(5)$ was improved to ${4\sqrt{10}}/{15}$ from $1$ in \cite{Catino2016}.
\endproof

\section{Some Lemmas}

In this section, we derive some inequalities which will play key roles
in proving our main results.
\begin{lem}\label{Lemma3-0}
Let $(M^n,g)$ be a gradient shrinking $\rho$-Einstein soliton satisfying
\eqref{Int2} with $n\geq4$. Then we have
\begin{align}\label{Lem-3-0-ineq}
0\geq&\Bigg\{\frac{n-2}{4(n-1)}Y(M,[g])-\frac{n-4}{2}\lambda
V(M)^{\frac{2}{n}}\nonumber\displaybreak[0]\\
&\quad-\sqrt{\frac{2(n-2)}{n-1}}\left(\int\limits_M\Big(|W|^2+\frac{8}{n(n-2)}|\mathring{{\rm
Ric}}|^2\Big)^{\frac{n}{4}}dv_g\right)^{\frac{2}{n}}\Bigg\}\Big(\int\limits_M|\mathring{{\rm
Ric}}|^\frac{2n}{n-2}dv_g\Big)^{\frac{n-2}{n}}\nonumber\displaybreak[0]\\
&\quad +\left(\frac{(n-4)^2}{4n(n-1)}-\frac{n-4}{2}\rho\right)\int\limits_MR|\mathring{{\rm
Ric}}|^2dv_g-\frac{(n-2)^2}{2n}\rho\int\limits_M|\nabla R|^2dv_g.
\end{align}
\end{lem}
\proof Using the inequality in Lemma \ref{weyl-est} and the Kato inequality $|\nabla\mathring{{\rm Ric}}|^2\geq|\nabla|\mathring{{\rm Ric}}||^2$ at
the point where $|\mathring{{\rm Ric}}|\neq0$, we have from
\eqref{2lemma-1} that
\begin{align*}
0\geq&-\frac{1}{2}\Delta_f|\mathring{{\rm Ric}}|^2+|\nabla|\mathring{{\rm
Ric}}||^2 +2\lambda|\mathring{{\rm
Ric}}|^2\\
&\quad-\sqrt{\frac{2(n-2)}{n-1}}\left(|W|^2+\frac{8}{n(n-2)}|\mathring{{\rm
Ric}}|^2\right)^{1/2}|\,\mathring{{\rm Ric}}|^2\\
&\quad +2\left(\rho-\frac{n-2}{n(n-1)}\right)R|\mathring{{\rm Ric}}|^2+(n-2)\rho
R_{,ij}\mathring{R}_{ij}.
\end{align*}
Integrating by parts over $M$ and using equation \eqref{lem1}, we obtain
\begin{align}\label{2Compact2}
0\geq&\int\limits_M\Big\{-\frac{1}{2}|\mathring{{\rm Ric}}|^2\Delta
f+|\nabla|\mathring{{\rm Ric}}||^2 +2\lambda|\mathring{{\rm
Ric}}|^2\nonumber\displaybreak[0]\\
&\qquad-\sqrt{\frac{2(n-2)}{n-1}}\left(|W|^2+\frac{8}{n(n-2)}|\mathring{{\rm
Ric}}|^2\right)^{1/2}|\,\mathring{{\rm Ric}}|^2\nonumber\displaybreak[0]\\
&\qquad +2\left(\rho-\frac{n-2}{n(n-1)}\right)R|\mathring{{\rm Ric}}|^2+(n-2)\rho
R_{,ij}\mathring{R}_{ij}\Big\}dv_g\nonumber\displaybreak[0]\\
=&\int\limits_M\Big\{|\nabla|\mathring{{\rm Ric}}||^2 -\frac{n
-4}{2}\lambda|\mathring{{\rm
Ric}}|^2-\sqrt{\frac{2(n-2)}{n-1}}\left(|W|^2+\frac{8}{n(n-2)}|\mathring{{\rm
Ric}}|^2\right)^{1/2}|\,\mathring{{\rm Ric}}|^2\nonumber\\
&\qquad+\left(\frac{n^2-5n+8}{2n(n-1)}-\frac{n-4}{2}\rho\right)R|\mathring{{\rm
Ric}}|^2+(n-2)\rho R_{,ij}\mathring{R}_{ij}\Big\}dv_g.
\end{align}
Notice that $\mathring{{\rm R}}_{ij,j}=\frac{n-2}{2n}R_{,i}$,  the last term in \eqref{2Compact2} is equal to
\begin{equation}\label{3-1-ineq1}
  (n-2)\rho\int\limits_M R_{,ij}\tilde{R}_{ij}dv_g~=~-\frac{(n-2)^2}{2n}\rho\int\limits_M|\nabla R|^2dv_g.
\end{equation}

As the Yamabe invariant $Y(M,[g])>0$ for any compact gradient shrinking $\rho$-Einstein soliton with $\rho<\frac 1{2(n-1)}$, we have the following Yamabe-Sobolev inequality by \eqref{Int3}
\begin{equation}\label{Int4}
\frac{n-2}{4(n-1)}Y(M,[g])\Big(\int_M|u|^{\frac{2n}{n-2}}\,dv_g\Big)^{\frac{n-2}{n}}\leq~\int_M|\nabla
u|^2\,dv_g+\frac{n-2}{4(n-1)}\int_M Ru^2\,dv_g.
\end{equation}
for any $u\in W^{1,2}(M)$. Replacing $u$ by  $|\mathring{{\rm Ric}}|$ in \eqref{Int4}, one
has
\begin{align*}
0\geq&~\frac{n-2}{4(n-1)}Y(M,[g])\Big(\int\limits_M|\mathring{{\rm
Ric}}|^{\frac{2n}{n-2}}dv_g\Big)^{\frac{n-2}{n}}-\frac{n
-4}{2}\lambda\int\limits_M|\mathring{{\rm
Ric}}|^2dv_g\\
&\quad -\sqrt{\frac{2(n-2)}{n-1}}\int\limits_M\left(|W|^2+\frac{8}{n(n-2)}|\mathring{{\rm
Ric}}|^2\right)^{1/2}\,|\mathring{{\rm Ric}}|^2dv_g\\
&\quad +\left(\frac{(n-4)^2}{4n(n-1)}-\frac{n-4}{2}\rho\right)\int\limits_MR|\mathring{{\rm
Ric}}|^2dv_g-\frac{(n-2)^2}{2n}\rho\int\limits_M|\nabla R|^2dv_g.
\end{align*}
As $\lambda>0$ and $n\geq 4$, the H\"{o}lder inequality then implies that
\begin{align*}
0\geq~&\Bigg\{\frac{n-2}{4(n-1)}Y(M,[g])-\frac{n-4}{2}\lambda
V(M)^{\frac{2}{n}}\\
&\quad -\sqrt{\frac{2(n-2)}{n-1}}\Big[\int\limits_M\Big(|W|^2+\frac{8}{n(n-2)}|\mathring{{\rm
Ric}}|^2\Big)^{\frac{n}{4}}dv_g\Big]^{\frac{2}{n}}\Bigg\}\Big(\int\limits_M|\mathring{{\rm
Ric}}|^\frac{2n}{n-2}dv_g\Big)^{\frac{n-2}{n}}\\
&\quad +\left(\frac{(n-4)^2}{4n(n-1)}-\frac{n-4}{2}\rho\right)\int\limits_MR|\mathring{{\rm
Ric}}|^2dv_g-\frac{(n-2)^2}{2n}\rho\int\limits_M|\nabla R|^2dv_g.
\end{align*}
\endproof

\begin{lem}\label{Lemma31}
Let $(M^n,g)$ be a gradient shrinking $\rho$-Einstein soliton satisfying
\eqref{Int2} with $n\geq4$. Then we have
\begin{align}\label{Lem-31}
0\geq&\Bigg\{\frac{n-2}{4(n-1)}Y(M,[g])-\frac{n-4}{2}\lambda
V(M)^{\frac{2}{n}}\nonumber\displaybreak[0]\\
&-\sqrt{\frac{2(n-2)}{n-1}}\left(\int\limits_M\Big(|W|^2+\frac{8}{n(n-2)}|\mathring{{\rm
Ric}}|^2\Big)^{\frac{n}{4}}dv_g\right)^{\frac{2}{n}}\Bigg\}\Big(\int\limits_M|\mathring{{\rm
Ric}}|^\frac{2n}{n-2}dv_g\Big)^{\frac{n-2}{n}}\nonumber\displaybreak[0]\\
&+\left(\frac{(n-4)^2}{4n(n-1)}-\frac{n-4}{2}\rho-\frac{(n-2)^2\rho}{n[1-2(n-1)\rho]}\right)\int\limits_MR|\mathring{{\rm
Ric}}|^2dv_g\nonumber\displaybreak[0]\\
&+\frac{(n-4)(n-2)^2\rho}{n^2[1-2(n-1)\rho]^2}\int\limits_M R\,{\rm
Ric}(\nabla f,\nabla f)dv_g.
\end{align}
\end{lem}
\proof
The proof is similar with Lemma \ref{Lemma3-0}. The only difference is that we estimate the last term in \eqref{2Compact2} as follows, by using \eqref{lem1}--\eqref{lem4}.
\begin{align*}
 &(n-2)\rho\int\limits_M R_{,ij}\mathring{R}_{ij}dv_g
=\frac{(n-2)^2}{2n}\rho\int\limits_M R\Delta
Rdv_g\nonumber\\
=&-\frac{(n-2)^2\rho}{n[1-2(n-1)\rho]}\int\limits_M
R\Big[-\frac{1}{2}\nabla R\nabla f+|\mathring{{\rm
Ric}}|^2+\frac{(1-n\rho)R-n\lambda}{n}R\Big]dv_g\nonumber\\
=&-\frac{(n-2)^2\rho}{n[1-2(n-1)\rho]}\int\limits_M
\Big[\frac{1}{4}R^2\Delta f+R|\mathring{{\rm
Ric}}|^2+\frac{(1-n\rho)R-n\lambda}{n}R^2\Big]dv_g\nonumber\displaybreak[0]\\
=&-\frac{(n-4)(n-2)^2\rho}{4n^2[1-2(n-1)\rho]}\int\limits_M
R^2\Delta fdv_g-\frac{(n-2)^2\rho}{n[1-2(n-1)\rho]}\int\limits_M
R|\mathring{{\rm Ric}}|^2dv_g\nonumber\\
=&\frac{(n-4)(n-2)^2\rho}{2n^2[1-2(n-1)\rho]}\int\limits_M R\nabla
R\nabla fdv_g-\frac{(n-2)^2\rho}{n[1-2(n-1)\rho]}\int\limits_M
R|\tilde{{\rm Ric}}|^2dv_g\nonumber\displaybreak[0]\\
=&\frac{(n-4)(n-2)^2\rho}{n^2[1-2(n-1)\rho]^2}\int\limits_M R\,{\rm
Ric}(\nabla f,\nabla
f)dv_g-\frac{(n-2)^2\rho}{n[1-2(n-1)\rho]}\int\limits_M R|\mathring{{\rm
Ric}}|^2dv_g.
\end{align*}
\endproof

Finally, we recall the following result which follows from the proof of \cite[Theorem 3.3]{Catino2016}, so we omit the proof.
\begin{lem}\label{Lemma32} Let $C(n)$ be defined in Lemma \ref{Lemma23}. Then for the
Einstein manifold $(M^n,g)$, $n\geq 4$, with positive scalar
curvature, we have
\begin{equation}\label{2Compact9}\aligned
&\left(\frac{n+1}{n-1}Y(M,[g])-\frac{8(n-1)}{n-2}C(n)\Big(\int\limits_M|W|^{\frac{n}{2}}dv_g\Big)^{\frac{2}{n}}\right)
\int\limits_M|\nabla|W||^2dv_g\\
&\qquad+\left(\frac{2}{n}Y(M,[g])-2C(n)\Big(\int\limits_M|W|^{\frac{n}{2}}dv_g\Big)^{\frac{2}{n}}\right)\int\limits_MR|W|^2dv_g\leq0.
\endaligned\end{equation}
\end{lem}

\section{Proof of main theorems}
\subsection{Proof of Theorem \ref{main-thm1}}

Since $\rho<0$, by Proposition \ref{Prop21} the scalar curvature $R>0$ on $(M^n,g)$.
Under the assumption of \eqref{main-ineq1}, in the case $n=5,6$, we get $\mathring{{\rm Ric}}=0$ from \eqref{Lem-3-0-ineq}
which says that $(M^n,g)$ is Einstein. When $n=4$, one can get $\nabla R=0$ from the last term in \eqref{Lem-3-0-ineq} and hence $R$ is constant. Furthermore, this shows that $(M^n,g)$ is Einstein for $n=4$ (see also \cite[Theorem 2.2]{HW2013}).
Thus, $\mathring{{\rm Ric}}=0$ holds for all $n=4,5,6$ and
the condition \eqref{main-ineq1} becomes
\begin{equation}\label{2Compact8}
\left(\int\limits_M|W|^{\frac{n}{2}}dv_g\right)^{\frac{2}{n}}
+\sqrt{\frac{(n-4)^{2}(n-1)}{8(n-2)}}\lambda V(M)^{\frac{2}{n}}\leq~\sqrt{\frac{n-2}{32(n-1)}}Y(M,[g]).
\end{equation}
Therefore $g$ is a Yamabe metric of the conformal class $[g]$, and
\begin{equation*}
  Y(M,[g])=V(M)^{\frac{2-n}n}\int\limits_MR~dv_g.
\end{equation*}
Integrating the equation \eqref{lem1} over $M$ gives that
\begin{equation*}
  \lambda V(M)^{\frac 2n}=\frac {1-n\rho}nV(M)^{\frac{2-n}n}\int\limits_MRdv_g=\frac {1-n\rho}nY(M,[g]).
\end{equation*}
Hence the pinching condition \eqref{2Compact8} implies that
\begin{align}\label{2Compact6}
\left(\int\limits_M|W|^{\frac{n}{2}}dv_g\right)^{\frac{2}{n}}
\leq &~\frac{8n-n^2-8+2n\rho(n-4)(n-1)}{4n\sqrt{2(n-1)(n-2)}}Y(M,[g])\nonumber\\
\leq &~\frac{8n-n^2-8}{4n\sqrt{2(n-1)(n-2)}}Y(M,[g])
\end{align}
as $\rho<0$. Substituting \eqref{2Compact6} into \eqref{2Compact9} and comparing the coefficients, we obtain
that the Weyl tensor $W=0$. Note that $(M^n,g)$ is Einstein and has positive scalar curvature. Therefore, we conclude that $(M^n,g)$ is isometric a quotient to the round sphere $\mathbb{S}^n$.

\subsection{Proof of Theorem \ref{1Th2}}
When $n\geq5$, since $M^n$ has positive Ricci curvature and $\rho>0$, we have
from\eqref{Lem-31}
\begin{align}\label{2Compact15}
0\geq&\Bigg\{\frac{n-2}{4(n-1)}Y(M,[g])-\frac{n-4}{2}\lambda
V(M)^{\frac{2}{n}}\nonumber\\
&-\sqrt{\frac{2(n-2)}{n-1}}\Big[\int\limits_M\Big(|W|^2+\frac{8}{n(n-2)}|\mathring{{\rm
Ric}}|^2\Big)^{\frac{n}{4}}dv_g\Big]^{\frac{2}{n}}\Bigg\}\Big(\int\limits_M|\mathring{{\rm
Ric}}|^\frac{2n}{n-2}dv_g\Big)^{\frac{n-2}{n}}\nonumber\displaybreak[0]\\
&+\Big[\frac{(n-4)^2}{4n(n-1)}-\frac{n-4}{2}\rho-\frac{(n-2)^2\rho}{n[1-2(n-1)\rho]}\Big]\int\limits_MR|\mathring{{\rm
Ric}}|^2dv_g.
\end{align}
The first condition \eqref{2th-1} implies that
\begin{equation}\label{rho-cond-1}
  \frac{(n-4)^2}{4n(n-1)}-\frac{n-4}{2}\rho-\frac{(n-2)^2\rho}{n(1-2(n-1)\rho)}\geq 0.
\end{equation}
Then combining \eqref{rho-cond-1} with the second condition \eqref{2th-2} of Theorem \ref{1Th2}, if at least one of them are strict inequality, we obtain from \eqref{2Compact15} that $(M^n,g)$ is Einstein. By Obata's \cite{Obata} theorem, $g$ is the Yamabe metric of the conformal class $[g]$.  Then, integrating the equation \eqref{lem1} gives
$$\lambda\,V(M)^{\frac{2}{n}}=\frac{1-n\rho}{n}V(M)^{\frac{2-n}{n}}\int\limits_M \,Rdv_g=\frac{1-n\rho}{n}Y(M,[g]).$$
Hence the condition \eqref{2th-2} becomes
\begin{equation}\label{2Compact20}
\Big(\int\limits_M|W|^{\frac{n}{2}}dv_g\Big)^{\frac{2}{n}}
\leq C_1(n,\rho)Y(M,[g]),
\end{equation}
where
\begin{equation}\label{def-C_1}
  C_1(n,\rho)=\frac{-n^2+8n-8}{4n\sqrt{2(n-1)(n-2)}}+\rho\sqrt{\frac{(n-4)^{2}(n-1)}{8(n-2)}}.
\end{equation}
Inserting \eqref{2Compact20} into the inequality \eqref{2Compact9}
yields
\begin{align}\label{2Compact21}
&\left(\frac{n+1}{n-1}-\frac{8(n-1)}{n-2}C_1(n,\rho)C(n)\right)\int\limits_M|\nabla|W||^2dv_g\nonumber\\
&\qquad+\left(\frac{2}{n}-2C_1(n,\rho)C(n)\right)\int\limits_MR|W|^2dv_g\leq0.
\end{align}
By a direct calculation, we can check that
\begin{equation}\label{2Compact22}\frac{n+1}{n-1}-\frac{8(n-1)}{n-2}
C_1(n,\rho)C(n)>0
\end{equation}
and
$$\frac{2}{n}-2C_1(n,\rho)C(n)>0.$$
hold for all $n\geq 5$ and $\rho$ satisfying the condition \eqref{2th-1}. Therefore \eqref{2Compact21} and the positivity of the scalar curvature imply that the Weyl tensor $W=0$. Then we conclude that $(M^n,g)$ is isometric to a quotient of a round sphere.

\subsection{Proof of Theorem \ref{1Th3}}

Assume that $(M^n,g)$ is a complete gradient shrinking $\rho$-Einstein soliton with $0<\rho<\frac 1{2(n-1)}$, bounded curvature, nonnegative sectional curvature and positive Ricci curvature. If the potential function $f$ is constant, then the equation \eqref{Int2} and the Bonnet-Myers theorem imply that the manifold is compact. If $f$ is not a constant, then reasoning as \cite[\S 3]{Catino2015}, there exists a hypersurface $\Sigma_0\subset M$ which is a regular level set of $f$ and the function $f$ only depends on the signed distance $r$ to $\Sigma_0$ on the whole manifold.  If the signed distance is bounded, then the potential function is also bounded (see Proposition 3.5 of \cite{Catino2015}). As the Bakry-Emery Ricci tensor $Ric+\nabla^2f$ is bounded from below by a positive constant, Theorem 1.4 of \cite{Wei-W} implies that the manifold is compact. If the signed distance $r$ is unbounded, Proposition 4.3 of \cite{Catino2015} implies that the scalar curvature $R$ is constant. Then our manifold $(M^n,g)$ is just a gradient shrinking Ricci soliton satisfying the assumption of \cite[Theorem 2]{Munt-Wang}. Hence, $(M^n,g)$ is compact.

\bibliographystyle{Plain}

\end{document}